\documentclass[11pt] {article} 
\usepackage[backref]{hyperref}
  \usepackage{amsmath}
    \usepackage{amssymb}
  \usepackage[dvipsnames]{xcolor}   
     \usepackage{pdfsync}

\newtheorem{theorem}{Theorem}[section]
\newtheorem{lemma}{Lemma}[section]

\newcommand{\eqnsection}{
   \renewcommand{\theequation}{\thesection.\arabic{equation}}
   \makeatletter
   \csname @addtoreset\endcsname{equation}{section} 
   \makeatother}

\def \ov{\overline}

\def \be{\begin{equation}}
\def \ee{\end{equation}}
\def \bt{\begin{theorem}} 
\def \et{\end{theorem}}
\def \bl{\begin{lemma}} 
\def \el{\end{lemma}}
\def \bea{\begin{eqnarray}}
\def \eea{\end{eqnarray}}
\def \bas{\begin{eqnarray*}}
\def \eas{\end{eqnarray*}}

\def \al{\alpha}
\def \bb{\beta}

\def \De{\Delta}
\def \ep{\epsilon}

\def \la{\lambda}

\def \om{\omega}
\def \Om{\Omega}

\def \vf{\varphi}
\def \si{\sigma}

\def \th{\theta}

\def \ze{\zeta}

\def \ff{\infty}
\def \wh{\widehat}
\def \wt{\widetilde}

\def \CC{{\cal C}}

\def \FF{{\cal F}}

\def \YY{{\cal Y}}

\def\b1{\mathbf 1}
\def \({\left(}
\def \){\right)}

\def \da{\downarrow }

\def \nn{\nonumber}
 
\def \Proof{\noindent{\bf Proof $\,$ }}

\def \bc{\begin{center} }
\def \ec{\end{center} }
\def \bs{\begin{slide} }
\def \es{\end{slide} }

\def\square{{\vcenter{\vbox{\hrule height.3pt
        \hbox{\vrule width.3pt height5pt \kern5pt
           \vrule width.3pt}
        \hrule height.3pt}}}}
\def\qed{{\hfill $\square$ \bigskip}}

\eqnsection

\begin{document}
 
 \title{ Ray-Knight theorems for the local times of rebirthed Markov processes }


 \author{   P.J. Fitzsimmons\,\,\,  Jay Rosen \thanks{Research of     Jay Rosen was partially supported by  grants from the Simons Foundation.   }}
\maketitle

  \abstract{We prove generalizations of the first and second Ray-Knight theorems, for a large class of non-symmetric strong Markov processes. These results link the local times of the Markov process with the squares of  associated Gaussian processes. This connection allows us to establish  results about the exact modulus of continuity (in the spatial variable) of the local times. Our approach is different from earlier treatments which were based on associated permanental processes rather than Gaussian processes. 
  
 The type of process with which we work can be described as follows. Start with a symmetric Markov process with finite lifetime; upon its death resurrect it at a place in the state space chosen at random, independent of the past. Continue in this way, resurrecting at each death, to obtain a recurrent process. The rebirthing procedure destroys the symmetry of the original process, leading to a large class of non-symmetric processes.
 
 The main results are illustrated by many examples.
 }

\section{Introduction}\label{sec-intro}

 Let $S$ a be locally compact space with a countable base.  
Let   $\YY=(\Om,  \FF_{t},  \YY_t,\th_{t},  P^x)$ be a transient symmetric Borel right process with state space $S$  and  continuous strictly positive  $p-$potential   densities  
\be
  u^{p}=\{ u^{p}(x,y), x,y\in S \},\quad p\geq 0,
\ee
 with respect to some $\si$--finite  positive measure $m$ on $S$.  
  Let $\ze=\inf\{t\,|\,\YY_t=\De \}$,  where $\De$ is the cemetery state for   $\YY$  and assume       that $\ze<\ff$   almost surely.

    We use the framework of rebirthed Markov processes given by Meyer \cite{Meyer}.  Let $\mu$ be a probability measure on $S$. Modify $\YY$ so that instead of going to $\De$ at the end of its lifetime it is immediately ``reborn'' with  probability $\mu $.  (I.e., the process goes to the set $B\subset S$ with probability $\mu(B) $, after which  it continues to evolve the way $\YY$ did, being reborn with probability measure  $\mu $ each time it dies.)   Denote this fully rebirthed process by $\wt Z\!=\!
(\wt \Om, \wt  \FF_{t}, \wt  Z_t, \wt \th_{t},  \wt P^y)$.   
 We show in    \cite[Theorem 7.1]{MRLIL}
that the process  
 $\wt   Z$ is a recurrent Borel right process with state space $ S$ and $p$--potential densities, 
\begin{equation}
w^{p}(x,y)=  u^{p}(x,y)+\(\frac{1}{p}-\int_S   u^{p}(x,z)\,dm(z)\)\frac{f(y)}{\|f\|_1},\qquad  p>0,\label{int.1}
\end{equation}
with respect to   $m$, where  
\begin{equation} \label{sint.2}
  f(y)=\int_S   u^{p}(x,y)\,d\mu(x),
\end{equation}
and the $L_1$ norm is taken with respect to $m$.
 Note that, even though $\YY$ is symmetric, $\wt   Z$ is not. 
 
Let  $\wt L_{t}^{y}$ denote the local time of $\wt Z$
  normalized so that,
\begin{equation}
  \wt E^{ x}\(\int_{0}^{\ff}  e^{-ps}\, d_{s}\wt L^{y}_s\)=   w^p(x,y). \label{int.3}
\end{equation}

Assume that $S\subseteq R^{1}$. In  \cite[Theorems  4.1, 4.2]{FMR}, we show that for a large class of $\YY$'s, for $d\in S$ and
 for some increasing function $\phi$,
   \be
   \limsup_{u\to 0}\frac{|\wt L_{t}^{d+u}-\wt L_{t}^{d}|}{\phi(x)  }=  2\(\wt L_{t}^{d}\)^{1/2},  \quad a.e.\,\,\, t,  \,\,\, \wt P^{y}\,\,
a.s.\label{irev.2}
\ee
for all $y\in S$, and  for closed intervals $\De\in S$  
\be
\lim_{h\to 0}\sup_{\stackrel{|u-v|\le h }{ u,v\in\De}} \frac{|\wt L_{t}^{u}-\wt L_{t}^{v}|}{\vf (u,v)  }= \sup_{u\in \De}  \(2\wt L_{t}^{u}\)^{1/2},  \quad a.e.\,\,\, t ,  \,\,\, \wt P^{y}\,\,
a.s.,\label{iurev.2c}
\ee
for all $y\in S$, for some   $\vf(u,v)$   such that $ \vf(u,v)\le  \wt\vf\(|u-v|\)$, where $\wt\vf$ is an increasing continuous function with $\wt\vf(0)=0$.
 
Note that  both these results only hold  for a.e. $t$. The goal of this paper is to examine such results when  $t$ is replaced by   random values,  $\wt T_{0}$
or $\wt \tau (t)$,  where 
\begin{equation}
\wt T_{0}=\inf\{s>0\,|\, \wt Z_{s}=0 \},\label{}
\end{equation}
and
\begin{equation}
\wt \tau (t)=\inf\{s>0\,|\,\wt L^{0}_{s}>t \}.\label{}
\end{equation} 
This is a non-trivial extension of previous work, requiring subtle conditional independence results (Lemmas \ref{lem-condind0}
and 
\ref{lem-ncondind}) to combine the re-birthing procedure with sampling at hitting times or at inverse local times.

We develop three different cases related to $\wt T_{0}$. In each of these we assume that the initial point $y\neq 0$  and that our renewal measure $\mu$ is supported away from $0$; that is, $0\notin  \overline{\text{supp}(\mu)}$.
  

 In Section \ref{sec-FRK} we prove Theorem \ref{theo-condind0}, a generalized First  Ray-Knight Theorem for  the local times $\wt L^{x}_{\wt T_{0}}$ of the non-symmetric Markov process $\wt Z$ (with values in $R^1$) which only involves squares of Gaussian processes.     We then use this to obtain exact moduli of continuity in $x$ for 
 $\wt L^{x}_{T_{0}}$ when $x$ is away from $0$.  

 In sub-Section \ref{subsec2} we assume that $S=R^{1}-\(0\)$ and set  
   \begin{equation}
 \wt T^{-}_{0}=\inf \{t\,|\, \wt Z_{t-}=0\}, \label{irkl.10}
  \end{equation}
where $\wt Z_{t-}$ is the left hand limit of $\wt Z_{t}$. We   show that for many examples of $\YY$ we have the analogs of 
(\ref{irev.2}) and (\ref{iurev.2c}) with $t $ replaced by $\wt T^{-}_{0}$, for $y \text{ and } d\neq 0$, and closed intervals $\De\in [0,1]\cap \{R^{1}-\(0\)\}$.

 In Section \ref{sec-atzero} we discuss the case where $\YY$ is an exponentially killed symmetric diffusion in $R^{1}$ which we then rebirth. (Note that since the rebirthed process is discontinuous at rebirth times,   it will no longer be a diffusion). We provide an exact local modulus for
 $\wt L^{x}_{\wt T_{0}}$ at $x=0$. Here, of course, $\wt L^{0}_{\wt T_{0}}=0$.

 
  We note that in \cite{Ray} Ray used the  first Ray-Knight theorem for Brownian motion (and more generally, diffusions) to obtain exact local and uniform moduli of continuity for the local times  $\{L_{T_{0}}^{x}, x\in R^{1}\}$.

  Finally in Section \ref{sec-NSRK} we prove Theorem \ref{theo-ITcond2}, a generalized Second Ray-Knight Theorem for    the local times  $\wt L^{x}_{\wt \tau (t)}$ of the non-symmetric Markov process $\wt Z$  which only involves squares of  Gaussian processes. 
We use this to exhibit many examples of $\YY$ for which  that the analogs of (\ref{irev.2}) and (\ref{iurev.2c}) hold with  $t $ replaced by $\wt \tau (t)$, $ d\neq 0$, and closed intervals $\De\in [0,1]\cap \{R^{1}-\(0\)\}$.


We will obtain our results by  generalizing the first and second Ray-Knight theorems to a large class of rebirthed Markov processes. 
 For any Markov process $X_{t}$ with local times $L_{t}^{x}$ we set
$ T_{0}=\inf\{s>0\,|\, X_{t}=0 \}$ and $ \tau (t)=\inf\{s>0\,|\, L^{0}_{s}>t \}$.
The  first and second Ray-Knight theorems  describe the laws of the Brownian local times $\{L_{T_{0}}^{x}, x\in R^{1}\}$ and $\{L_{\tau(t)}^{x}, x\in R^{1}\}$.
      \, These theorems can be formulated in terms of a connection between  $\{L_{T_{0}}^{x}, x\in R^{1}\}$, $\{L_{\tau(t)}^{x}, x\in R^{1}\}$ and squares of Gaussian processes,  \cite[Theorems 2.6.3 and 2.7.1]{book}. 
In \cite{EKMRS} such connections were generalized to all strongly symmetric Markov processes with finite 1-potential densities, see \cite[Sections 8.1.1 and 8.2]{book}.  As mentioned already, the rebirthed Markov processes we consider in this paper are {\em not} symmetric.

 In \cite[Corollary 3.5]{EK} the authors  provide connections for any recurrent Markov process between $  L^{x}_{  \tau (t)}$ and  permanental processes. 
We have not been able to obtain exact uniform moduli of continuity using 
 permanental processes. The innovation of the  approach taken in this paper  is that we connect $  L^{x}_{ T_{0}}$ and  $  L^{x}_{  \tau (t)}$ with squares of  Gaussian processes and use this to obtain exact local and uniform moduli of continuity for those  local times.




\section{A Generalized First Ray-Knight Theorem for rebirthed Markov Processes}\label{sec-FRK}

As described in the Introduction, let $S$ a be locally compact space with a countable base.  
Let   $\YY=(\Om,  \FF_{t},  \YY_t,\th_{t},  P^x)$ be a transient symmetric Borel right process with state space $S$  and  continuous strictly positive  $p-$potential   densities  
\be
  u^{p}=\{ u^{p}(x,y), x,y\in S \},\quad p\geq 0,
\ee
 with respect to some $\si$--finite  positive measure $m$ on $S$. Let $\ze=\inf\{t\,|\,\YY_t=\De \}$,  where $\De$ is the cemetery state for   $\YY$  and assume       that $\ze<\ff$   almost surely.       It follows from  \cite[Lemma 3.3.3]{book} that $ \{ u^{p}(x,y), x,y\in S \} $   is positive definite and therefore is  the covariance  of a Gaussian process.

 \medskip Let $\mu$ be a probability measure on $S$. Following \cite{Meyer} we       modify $\YY$ so that instead of going to $\De$ at the end of its lifetime it is immediately reborn with  measure $\mu $.  (I.e., the process goes to the set $B\subset S$ with probability $\mu(B) $, after which  it continues to evolve the way $\YY$ did, being reborn with probability   $\mu $ each time it dies.)  
  
   We denote this rebirthed process by $\wt Z\!=\!
(\wt \Om, \wt  \FF_{t}, \wt  Z_t, \wt \th_{t},\wt P^x)$.   
We have $\wt \Om=\Om^{ \mathbf N}$ with elements $\wt\om=( \om_{1},\om_{2}, \ldots)$
 and $\wt P^{y}=  P_{1}^{y}\times_{n=2}^{\ff}P_n^{\mu}$ on $( \Om,   \FF  )^{ \mathbf N}$, where   for any probability measure $P$, $\{P_i,i\ge 1 \}$ are independent copies of $P$.

We show in    \cite[Theorem 7.1]{MRLIL}
that the process  
 $\wt   Z$ is a recurrent Borel right process with state space $ S$ and $p$--potential densities, 
 given by 
(\ref{int.1}) and (\ref{sint.2}).
It should be noted that
 \begin{equation} \label{mar.1}
  f(y)\le u^{p}(y,y)<\infty,\qquad \forall y\in S,
\end{equation}  
because for symmetric processes,  $u^{p}(x,y)\leq u^{p}(y,y)$,
  $\forall x, y\in S $, and $\mu(S)=1$.
Note also  that $w^{p}(x,y)$, and hence $\wt Z$, is not symmetric.

 Let  $\{\wt L_{t}^{y}, y\in S,t\in R_+ \}$ denote the local time of $\wt Z$
  normalized so that,
\begin{equation}
  \wt E^{ x}\(\int_{0}^{\ff}   e^{-ps}\, d_{s}\wt L^{y}_s\)=   w^p(x,y), \label{az.1vw}
\end{equation}
 and let $0$ be  a fixed point in $S$.

  In the following we take our initial point $y\neq 0$ and our renewal measure $\mu$ supported away from $0$, that is, $0\notin  \overline{\text{supp}(\mu)}$.
  
   By our assumption that the $p$-potential   densities  for $\YY$ are  continuous and strictly positive, it follows from \cite[(3.107)]{book} that  $0$ is not polar for $\YY$, that is  $P^{x}(T_{0}<\ff)>0$ for all $x\in S$, where $T_{0}=\inf\{s>0\,|\, \YY_{s}=0 \}$. Thus it is not polar for $\wt Z$, and since $\wt Z$ is recurrent, it follows from \cite[(10.39)]{S} that $\wt P ^{y}\( \wt T_{0}(\wt\om)<\ff  \)=1$, where   $\wt T_{0}=\inf\{s>0\,|\, \wt Z_{s}=0 \}$. 



  Let  $L_t^x=\{L^{x}_{t},y\in S,t\in R^+\}  $ be the local time for $\YY$ normalized so that
  \begin{equation}\label{2.2}
  E^{z}\( L_{\ff}^{x}\)= u^{0}(z,x), \quad \forall z,x\in S.
  \end{equation}
  We require that $L^{x}_{t}  $ is jointly continuous.
Let $    \wt L^{x}_{t}$ be the local time for  $\wt Z$ 
 and let 
 $\ze_{n}=\ze_{n}(\wt\om)=\sum_{j=1}^{n}\ze (\om_{j})$ and $\ze_{0}=0$.
We show in \cite[Lemma 2.1]{FMR}
that for each $r\ge 1$,
  \begin{equation}
 \wt L^{x}_{t}(\wt\om)=\sum^{r-1}_{i=1}L^{x}_{\ze (\om_{i})}(\om_{i})+ L^{x}_{t-\ze_{r-1}}(\om_{r}), \quad \forall t\in(\ze_{r-1},\ze_{r}], 
 \quad a.s.\,,\label{74.1}
 \end{equation} 
 so that if $\ze_{r-1}<\wt T_{0}=\wt T_{0}(\wt\om)<\ze_{r}$ then
   \begin{equation}
 \wt L^{x}_{\wt T_{0}}(\wt\om)=\sum^{r-1}_{i=1}  L^{x}_{\ze(\om_{i})}(\om_{i}) +     L^{x}_{T_{0}(\om_{r})}(\om_{r}).\label{74.2}
 \end{equation}
We note that we cannot have $\wt T_{0}(\wt\om)=\ze_{r}$ for any $r$ since $\wt Z_{\wt T_{0}(\wt\om)}=0$ while $\wt Z_{\ze_{r}}\in \text{supp}(\mu)$ and by assumption $0\notin \overline {\text{supp}(\mu)}$.

 For each $i$ let
 \begin{equation}
 \mathcal{L}_{i}=\{L^{x}_{\ze(\om_{i})}(\om_{i}), x\in S\}\label{}
 \end{equation}
 and
  \begin{equation}
 \mathcal{L}'_{r}=\{L^{x}_{T_{0}(\om_{r})}(\om_{r}), x\in S\}.\label{}
 \end{equation}

 \bl\label{lem-condind0}
  For each $r\geq 2$, and all  Borel sets $B_1, B_2,\ldots, B_r$  in $C(S, R^1)$,
  \begin{eqnarray}
 && \wt P ^{y}\( \mathcal{L}_{1}\in B_1, \ldots,  \mathcal{L}_{r-1}\in B_{r-1},  \mathcal{L}'_{r} \in B_r \,\Big |\,    \ze_{r-1}<\wt T_{0}< \ze_{r} \)
 \label{74.5mp}
 \\
 &&= P_{1} ^{y}\( \mathcal{L}_{1}\in B_1 \,\Big |\,   \ze(\om_{1})<T_{0}(\om_{1})  \) \prod^{r-1}_{i=2}    P_{i} ^{\mu}\(\mathcal{L}_{i}\in B_i\,\Big |\,  \ze(\om_{i})<T_{0}(\om_{i}) \)\nn\\
 &&\hspace{2 in}  \times   P_{r} ^{\mu}\(  \mathcal{L}'_{r}\in B_r\,\Big |\,  T_{0}(\om_{r})< \ze (\om_{r})\).
\nn
 \end{eqnarray}
 \el

\textbf{Proof: } Let $\mathcal{C}$ be a countable subset of $S$ with compact closure.    For any function $f(x)$ and $i=1,\ldots, r$,  let $(\nu_{i},f)=\sum_{x\in \mathcal{C}}a_{i, x} f(x)$ where the $a_{i, x}\in R^{1}$ are arbitrary except that we assume $\sum_{x\in \mathcal{C}}|a_{i, x}|<\ff$ for each $i=1,\ldots, r$. By (\ref{74.2})
 \begin{eqnarray}
 &&\wt E ^{y}\(e^{-\sum^{r-1}_{i=1} (\nu_{i},   L^{\cdot}_{\ze(\om_{i})}(\om_{i}))-  (\nu_{r},  L^{\cdot}_{T_{0}(\om_{r})}(\om_{r}))}   ;    \ze_{r-1}<\wt T_{0}<\ze_{r} \)
\nn
 \\
  &&=\wt E^{y}\(   e^{-\sum^{r-1}_{i=1} (\nu_{i},  L^{\cdot}_{\ze(\om_{i})}(\om_{i}))}  ;    \ze_{r-1}<\wt T_{0}   \)       E^{\mu}\(  e^{-  (\nu_{r},  L^{\cdot}_{T_{0}(\om_{r})}(\om_{r}))}; T_{0}(\om_{r})<\ze (\om_{r})  \)
 \nonumber\\
  &&=\wt E^{y}\( \prod^{r-1}_{i=1} e^{-  (\nu_{i},  L^{\cdot}_{\ze(\om_{i})}(\om_{i}))}1_{\{\ze (\om_{i})<T_{0}(\om_{i})\}}     \)         E^{\mu}\(  e^{-  (\nu_{r},  L^{\cdot}_{T_{0}(\om_{r})}(\om_{r}))}; T_{0}(\om_{r})< \ze (\om_{r})  \)\nn\\
 &&=  E_{1}^{y}\(  e^{-  (\nu_{1},  L^{\cdot}_{\ze(\om_{1})}(\om_{1}))} ;  \ze(\om_{1})<T_{0}(\om_{1})   \) \prod^{r-1}_{i=2}   E_{i}^{\mu}\( e^{-  (\nu_{i},  L^{\cdot}_{\ze(\om_{i})}(\om_{i}))} ; \ze(\om_{i})<T_{0}(\om_{i})\) \nonumber\\
 &&\hspace{1 in} \times    E_{r}^{\mu}\(  e^{-  (\nu_{r},  L^{\cdot}_{T_{0}(\om_{r})}(\om_{r}))}; T_{0}(\om_{r})< \ze (\om_{r})  \).
 \label{74.3}
 \end{eqnarray}
 
Since, on setting all $\nu_{i}=0$
 \bea
&&\wt P ^{y} \(  \ze_{r-1}<\wt T_{0}< \ze_{r}\)\nn\\
&&=   P_{1} ^{y} \( \ze(\om_{1})<T_{0}(\om_{1}) \)\prod^{r-1}_{i=2}   P_{i} ^{\mu} \( \ze(\om_{i})<T_{0}(\om_{i})\)    P_{r} ^{\mu} \(T_{0}(\om_{r})< \ze (\om_{r})\),\nn
 \eea
 we have that
 \begin{eqnarray}
  &&\wt E ^{y}\(e^{-\sum^{r-1}_{i=1} (\nu_{i},   L^{\cdot}_{\ze(\om_{i})}(\om_{i}))-  (\nu_{r},  L^{\cdot}_{T_{0}(\om_{r})}(\om_{r}))}    \,\Big |\,    \ze_{r-1}<\wt T_{0}< \ze_{r} \)
\nn
 \\
 &&= E_{1}^{y}\(  e^{-  (\nu_{1},  L^{\cdot}_{\ze(\om_{1})}(\om_{1}))} \,\Big |\,  \ze(\om_{1})<T_{0}(\om_{1})   \)  \label{74.3n}\\
 &&\hspace{.8 in} \times\prod^{r-1}_{i=2}    E_{i}^{\mu}\( e^{-  (\nu_{i},  L^{\cdot}_{\ze(\om_{i})}(\om_{i}))} \,\Big |\, \ze(\om_{i})<T_{0}(\om_{i})\) \nonumber\\
 &&\hspace{1 in} \times     E_{r}^{\mu}\(  e^{-  (\nu_{r},  L^{\cdot}_{T_{0}(\om_{r})}(\om_{r}))}\,\Big |\, T_{0}(\om_{r})<\ze (\om_{r})  \).
 \nonumber
 \end{eqnarray}
 
By Laplace inversion this gives   (\ref{74.5mp}) with the $B_{i}$ restricted to $\mathcal{C}$, and the general case then follows by continuity.
 \qed

  Let 
 $L^{x}_{i, t}$ for $i=1, 2, \ldots$ be independent copies of $L^{x}_{t}$.
  For any set $C$, let $F(C)$ denote the set of real--valued functions
$f$ on $C$. Define the evaluations
$\it{i}_{ x}:F(C)\mapsto R^{{\rm 1}}$ by $\it{i}_{ x}( f)=f( x)$. We use
$\mathcal{M}( F(C))$ to denote  the smallest $\si$-algebra for which the
evaluations
$\it{i}_{ x}$ are Borel measurable for all
$x\in C$. ($\mathcal{M}( F(C))$ is generally referred to  as the
$\si$-algebra  of cylinder sets in $F(C)$.) 
As in \cite[Theorem 2.1]{FMR}, using   Lemma \ref{lem-condind0}  we obtain the following Theorem.

 \bt \label{cor-condind0} If for some $r\geq 2$ and measurable set of functions $B\in \mathcal{M}( F(S))$
 \begin{equation}
\(P_{1}^{y}\times \prod^{r}_{i=2} P_{i}^{\mu} \)  \( \sum^{r-1}_{i=1} L^{\cdot}_{i, \ff}+ L^{\cdot}_{r, T_{0}}\in B\)=1,  \label{74.6}
 \end{equation}
 then  
  \begin{equation}
\wt P^{y} \(\wt L^{\cdot}_{\wt T_{0}}(\wt\om)\in B\,\Big |\,    \ze_{r-1}<\wt T_{0}< \ze_{r}\)=1.\label{74.7}
 \end{equation}
 \et
 
  Note that on the right hand side of (\ref{74.5mp}), $L^{\cdot}_{r, T_{0}}$ was conditioned on $ T_{0}(\om_{r})< \ze (\om_{r})$,      so that the process actually reached $0$. Now that we removed the conditioning it is possible that $ T_{0}(\om_{r})=\ff$, so that we should think of $L^{\cdot}_{r, T_{0}}$ as the total local time of the process $\wt \YY$ obtained by killing $ \YY$ the first time it hits $0$. If $u^{0}(x,y)$ is the potential density      for  $ \YY$,    then by \cite[(4.165)]{book}  the potential density for the process $\wt \YY$  is
 \begin{equation}
 \wt   u^{0}(x,y)=u^{0}(x,y)- \frac{u^{0}(x,0)u^{0}(0,y)}{u^{0}(0,0)}.\label{74.new}
 \end{equation} 
 
 
 Let $\eta_{i, 0} (x)$, i=1,2,\ldots be independent  Gaussian processes with covariance $ u^{0}(x,y)$ and let  $\wt \eta (x)$ be a Gaussian process with covariance $ \wt   u^{0}(x,y)$ independent of the $\eta_{i, 0} (x)$.
 
The following is our generalization of the  first Ray-Knight theorem to rebirthed Markov processes.

 \bt\label{theo-condind0} For all $r\geq 1$ set
\be \ov G_{r, s}\(x\)=\sum_{i=1}^{r-1} \frac{1}{2}(\eta_{i, 0 }(x)+s)^2+\frac{1}{2}(\wt \eta (x)+s)^2, \label{2.21phh}
\ee
and $ P_{\ov G_{r, s}} =  \prod_{i=1}^{r-1}  P_{\eta_{i, 0} } \times   P_{\wt \eta} $. If
\begin{equation}
P_{\ov G_{r, s}}\(\ov G_{r, s}(\cdot)\in B\)=1,  \label{70.6nb}
 \end{equation}
  for a measurable set of functions $B\in \mathcal{M}( F(S))$, then 
  \begin{equation}
\(\wt P^{y}\times P_{\ov G_{r, s}}\) \(\wt L^{\cdot}_{\wt T_{0}}(\wt\om)+\ov G_{r, s}(\cdot)\in B\,\Big |\,    \ze_{r-1}<\wt T_{0}\leq \ze_{r}\)=1.\label{70.7n}
 \end{equation}
\et

 \Proof  We use the Eisenbaum Isomorphism Theorem, \cite[Theorem 8.1.1]{book}, which gives the relationships,   
\bea\lefteqn{
\Big\{ L^x_{i,\ff}+\frac{1}{2}(\eta_{0,i} (x)+s)^2\,\,;\,x\in S\,,\,P_{i}^y\times
P_{\eta_{0,i}}\Big\}\label{4it1.2}}\\ &&\qquad\stackrel{law}{=}
\Big\{\frac{1}{2}(\eta_{0,i} (x)+s)^2\,\,;\,x\in S\,,\,\(1+\frac{\eta_{0,i} (y)}s\)P_{\eta_{0,i}}\Big\},\nn
\eea
where $P_{\eta_{0, i}}$ is the probability for $\eta_{0, i}$ and $P_{i}^y$ the probability for $\YY$ started at $y$, 
and, as explained above (\ref{74.new}),
\bea\lefteqn{
\Big\{  L^x_{r, T_{0}}+\frac{1}{2}(\wt \eta (x)+s)^2\,\,;\,x\in S\,,\,  P_{r}^y\times
P_{\wt \eta}\Big\}\label{4it1.2a}}\\ &&\qquad\stackrel{law}{=}
\Big\{\frac{1}{2}(\wt \eta (x)+s)^2\,\,;\,x\in S\,,\,(1+\frac{\wt \eta (y)}s)P_{\wt \eta}\Big\}\nn
\eea
where $P_{\wt \eta}$ is the probability for $\wt \eta$.

Set $\eta_{0,i} (\mu)=\int \eta_{0,i} (x)\,d\mu(x)$,  $\wt \eta (\mu)=\int \wt \eta (x)\,d\mu(x)$.   It follows from the preceding two relationships     that  for each  $r\geq 2$,
\bea
&&\label{4it1.3}
\Big\{ \sum^{r-1}_{i=1} L^{x}_{i, \ff}+  L^{x}_{r, T_{0}}+\sum_{i=1}^{r-1}\frac{1}{2}(\eta_{0,i}(x)+s)^2+\frac{1}{2}(\wt \eta(x)+s)^2\,;\,\, x\in S,\,\,  \\
&&\hspace{1.5 in} \( P^{y} \times\prod_{i=2}^{r}P_{i}^{\mu}  \)     \times \prod_{i=1}^{r-1}  P_{\eta_{0,i} } \times   P_{\wt \eta } \Big\}\nn\\
 &&\stackrel{law}{=}
\Big\{ \sum_{i=1}^{r-1}\frac{1}{2}(\eta_{0,i}(x)+s)^2+\frac{1}{2}(\wt \eta(x)+s)^2\,\,\,\,\,;\,x\in S\,,\nn\\
&& \hspace{.5 in}\qquad     \,(1+\frac{\eta_{ 0,1}(y)}s)P_{\eta_{0,1} } \times\prod_{i=2}^{r-1}  (1+\frac{\eta_{0,i}(\mu)}s)P_{\eta_{0,i} }\times   ( (1+\frac{\wt \eta(\mu)}s)P_{\wt \eta }) \Big\}.\nn
\eea
Theorem \ref{theo-condind0} for  $r=1$ is precisely (\ref{4it1.2a}).


 For $r>1$ it follows from (\ref{70.6nb}) and (\ref{4it1.3}) that 
 \be
 \(P_{1}^{y} \times\prod_{i=2}^{r}P_i^{\mu}  \times P_{\ov G_{r, s}}\)
\( \sum^{r-1}_{i=1} L^{\cdot}_{i, \ff}+ L^{\cdot}_{r, T_{0}}+\ov  G_{r,s}(\cdot)\in B\)=1. \label{}
\ee
Thus for $P_{\ov G_{r, s}}$ almost every $\om'$
 \be
 \(P_{1}^{y} \times\prod_{i=2}^{r}P_i^{\mu}   \)
\( \sum^{r-1}_{i=1} L^{\cdot}_{i, \ff}+ L^{\cdot}_{r, T_{0}}+\ov  G_{r,s}(\cdot, \om')\in B\)=1. \label{}
\ee
 
  Theorem \ref{theo-condind0} then follows    from Theorem \ref{cor-condind0}.\qed

 As in \cite[Section 4]{FMR} this leads to
  \bt\label{theo-localmod0} Assume that $S\subseteq R^{1} $.  Let $\eta_{0} (t)$ be a Gaussian process with covariance $ u^{0}(x,y)$ and let $\wt  \eta  (t)$ be a Gaussian process with covariance $\wt  u^{0}(x,y)$. Set,
  \begin{equation} \label{}
  \si_{0}^{2} (x,y)=  u^{0}(x,x)+ u^{0}(y,y)-2 u^{0}(x,y) 
\end{equation}
and,
  \begin{equation} \label{}
 \wt   \si_{0} ^{2} (x,y)= \wt  u^{0}(x,x)+ \wt  u^{0}(y,y)-2  \wt  u^{0}(x,y). 
\end{equation} 
If,  for some $y_{0}\neq 0$
   \begin{equation}
 \lim_{x,y\da y_{0}} \frac{  \wt   \si_{0} ^2(x,y)}{\si_{0} ^2(x,y)}=1\label{41.4am},
   \end{equation}
  and   $\{\eta_{0} 
  (t)\}$ has a local modulus of continuity,  
  \begin{equation} \label{41.7m}
   \limsup_{x\to 0} \frac{|\eta_{0} (y_{0}+x)-\eta_{0}  (y_{0})|} {   \phi(x)} =  1    \qquad a.s.,
\end{equation}
for some increasing function $\phi$,  then  
   \be
   \limsup_{x\to 0}\frac{|\wt L_{\wt T_{0}}^{y_{0}+x}-\wt L_{\wt T_{0}}^{y_{0}}|}{\phi(x)  }= \(2\wt L_{\wt T_{0}}^{y_{0}}\)^{1/2},  \quad   \,\,\, \wt P^{y}\,\,
a.s.\label{4rev.2}
\ee
for all $y\in S$, $y\neq 0$.
\et

In \cite[Section 5]{FMR} we give several examples for L\'{e}vy processes and diffusions.  Let $Y=\{Y_{t};t\ge 0 \}$ be a real valued symmetric 
L\'evy process.
In Case 1 we let  $\YY=Y_{1}$ be the process with state space $S=R^{1}$, that is obtained by killing  $Y$ at the end of an independent exponential time with mean $1/\bb$. In this case  our $u^{0}(x,y)$ corresponds to the continuous $\bb>0$ potential density $u^{\bb}(x-y)$ of the L\'{e}vy process $Y$. 

Also, $\wt  u^{0}(x,y)$ corresponds to the continuous  $\bb>0$ potential density  $v^{\bb}(x,y)$ of Case 2. The fact that (\ref{41.4am}) holds is proven in \cite[(5.31)]{FMR},  and in that section we show that (\ref{41.7m}) holds with $ \phi(x)$ asymptotic to $2(\si_{0} )^{2} (x)\log\log 1/ |x|$ at $0$.

Cases 2 and 3 do not apply  since in those cases the processes  $\YY=Y_{2},Y_{3}$  have state space $S=R^{1}-\{0\}$,
which does not include $0$. 

Similar remarks apply to our examples for diffusions. Case 4 applies, but not Cases 5 and 6 which have state space $(0,\ff)$.

 Also, as in \cite[Section 4]{FMR}, there is a similar result for the exact uniform modulus of continuity.

\bt\label{theo-unifmodi}    
Assume that $S\subseteq R^{1} $,
  \begin{equation} \label{i.12}
 \lim_{x\to 0}\sup_{|u-v|\le x}(\si_{0} )^{2} (u,v)\log 1/ |u-v|=0, 
\end{equation}
and similarly for $(\wt \si_{0} )^{2}$.

  Let $\vf(u,v)$ be   such that for some constant $0<C<\ff$,
\begin{equation} \label{i.13}
  \vf(u,v)\le  \wt\vf\(C|u-v|\),
\end{equation}
 where $\wt\vf$ is an increasing continuous function with $\wt\vf(0)=0$.     If
 \be
 \lim_{h\to 0}\sup_{\stackrel{|u-v|\le h }{ u,v\in\De}} \frac{ ( \eta _{0 }(u)-\eta _{0}(v)) a_{0} + (\wt\eta (u)-\wt\eta (v)) a_{1} }{    \    \vf (u,v) } =1 \label{i.11} ,
\ee
 for all   intervals $\De$ in $[0,1]\cap \(S-\{0\}\)$ and all     $\{a_0,a_1\}$  with  $a_0^2+a_1^2=1$     then,  
 
  \be
\lim_{h\to 0}\sup_{\stackrel{|u-v|\le h }{ u,v\in\De}} \frac{|\wt L_{\wt T_{0}}^{u}-\wt L_{\wt T_{0}}^{v}|}{\vf (u,v)  }= \sup_{u\in \De}  \(2\wt L_{\wt T_{0}}^{u}\)^{1/2},  \quad   \,\,\, \wt P^{y}\,\,
a.s.,\label{i.14}
\ee
for all $y\in S$, $y\neq 0$.  
\et 

 Using \cite[Section 5]{FMR} we can see that this Theorem holds in Case1 with
$\vf (u,v)$ asymptotic to $2(\si_{0} )^{2} (u,v)\log 1/ |u-v|$.

\subsection{Case 2: A Generalized First Ray-Knight Theorem for rebirthed Markov Processes which do not hit $0$}\label{subsec2}

As mentioned before, in Case 2 the process  $\YY=Y_{2}$  has state space $S=R^{1}-\{0\}$,
which does not include $0$, so the discussion so far does not apply. However, we can still obtain interesting results. Here, of course, we will be starting at $y\neq 0$  and our renewal measure $\mu$ is supported away from $0$, that is $0\notin \overline {\text{supp}(\mu)}$. Since $\YY$  has state space $S=R^{1}-\{0\}$, $\wt Z$ will never hit $0$.   For any function  $\{f_{s}, s\geq 0\}$ with left hand limits,  let $f_{t^{-}}$ denote the left limit of $f$ at $t$.    We set
  \begin{equation}
\wt   T^{-}_{0}( \wt  \om)=\inf \{t\,|\, \wt Z_{t-}(\wt   \om)=0\}, \label{rkl.10}
  \end{equation}
  and consider $\wt  L^{x}_{\wt  T^{-}_{0}(\wt \om)}(\wt \om)$.


Using (\ref{74.1}), as in (\ref{74.2})
we see that if $\ze_{r-1}<\wt  T^{-}_{0}=\wt  T^{-}_{0}(\wt\om)\leq \ze_{r}$ then
   \begin{equation}
 \wt L^{x}_{\wt  T^{-}_{0}}(\wt\om)=\sum^{r-1}_{i=1}  L^{x}_{\ze(\om_{i})}(\om_{i}) +     L^{x}_{T^{-}_{0}(\om_{r})}(\om_{r}).\label{74.2r}
 \end{equation}

 For each $i$ let
 \begin{equation}
 \mathcal{L}_{i}=\{L^{x}_{\ze(\om_{i})}(\om_{i}), x\in S\}\label{}
 \end{equation}
 and
  \begin{equation}
 \mathcal{L}'_{r}=\{L^{x}_{T^{-}_{0}(\om_{r})}(\om_{r}), x\in S\}.\label{}
 \end{equation}

 \bl\label{lem-condind0-}
  For each $r\geq 2$, and all  Borel sets $B_1, B_2,\ldots, B_r$  in $C(S, R^1)$,
  \begin{eqnarray}
 && \wt P ^{y}\( \mathcal{L}_{1}\in B_1, \ldots,  \mathcal{L}_{r-1}\in B_{r-1},  \mathcal{L}'_{r} \in B_r \,\Big |\,    \ze_{r-1}<\wt  T^{-}_{0}\leq \ze_{r} \)
\nn
 \\
 &&= P_{1} ^{y}\( \mathcal{L}_{1}\in B_1 \,\Big |\,   \ze(\om_{1})<T^{-}_{0}(\om_{1})  \) \prod^{r-1}_{i=2}    P_{i} ^{\mu}\(\mathcal{L}_{i}\in B_i\,\Big |\,  \ze(\om_{i})<T^{-}_{0}(\om_{i}) \)\nn\\
 &&\hspace{1.6 in}  \times   P_{r} ^{\mu}\(  \mathcal{L}'_{r}\in B_r\,\Big |\,  T^{-}_{0}(\om_{r})\leq  \ze (\om_{r})\).
 \label{rkl.13}
 \end{eqnarray}
 \el
The proof of this Lemma follows as before. 

 Let 
 $L^{x}_{i, t}$ for $i=1, 2, \ldots$ be independent copies of $L^{x}_{t}$.
As in \cite[Theorem 2.1]{FMR}, using   Lemma \ref{lem-condind0}  we obtain the following Theorem.

 \bt \label{cor-condind0-} If for some $r\geq 2$ and measurable set of functions $B\in \mathcal{M}( F(S))$
 \begin{equation}
\(P_{1}^{y}\times \prod^{r}_{i=2} P_{i}^{\mu} \)  \( \sum^{r-1}_{i=1} L^{\cdot}_{i, \ff}+ L^{\cdot}_{r, T^{-}_{0}}\in B\)=1,  \label{rkl.14}
 \end{equation}
 then  
  \begin{equation}
\wt P^{y} \(\wt L^{\cdot}_{\wt  T^{-}_{0}}(\wt\om)\in B\,\Big |\,    \ze_{r-1}<\wt  T^{-}_{0}\leq  \ze_{r}\)=1.\label{rkl.15}
 \end{equation}
 \et
 
 Recall that in  Case 1 we let  $\YY=Y_{1}$ be the process   obtained by killing the L\'{e}vy process  $Y$ at the end of an independent exponential time. For ease of notation we set $\bar Y=Y_{1}$. Let $\bar T_{0}=\inf\{s>0\,|\, \bar Y_{s}=0 \}$ and $\bar T^{-}_{0} =\inf \{s>0\,|\, \bar Y^{-}_{s}=0\}$. It follows from the quasi-left continuity of the L\'{e}vy process  $Y$, \cite[Chapter 1, Proposition 7]{Bertoin}, that if 
$\bar T^{-}_{0}<\ze$ then $\bar T^{-}_{0}=\bar T_{0}$. 

Since   $\YY=Y_{2}$ in Case 2  is obtained by killing  $\bar Y$ the first time it hits $0$, it follows that 
\begin{equation}
\{ T^{-}_{0}(\om_{r})\leq  \ze (\om_{r})\}=\{ T^{-}_{0}(\om_{r})=\ze (\om_{r})\}.\label{2.qc}
\end{equation}

 Note that on the right hand side of (\ref{rkl.13}), $L^{\cdot}_{r, T^{-}_{0}}$ was conditioned on $ T^{-}_{0}(\om_{r})\leq \ze (\om_{r})$, or equivalently, by  (\ref{2.qc}),  
 $ T^{-}_{0}(\om_{r})= \ze (\om_{r})$,  so that the process actually converged to $0$. Now that we removed the conditioning it is possible that $ T^{-}_{0}(\om_{r})=\ff$, so that in either case  $L^{\cdot}_{r, T^{-}_{0}}=L^{\cdot}_{r, \ff}$. Thus our Lemma says that  
if for some $r\geq 2$ and measurable set of functions $B\in \mathcal{M}( F(S))$
 \begin{equation}
\(P_{1}^{y}\times \prod^{r}_{i=2} P_{i}^{\mu} \)  \( \sum^{r}_{i=1} L^{\cdot}_{i, \ff}\in B\)=1,  \label{rkl.16}
 \end{equation}
 then  
  \begin{equation}
\wt P^{y} \(\wt L^{\cdot}_{\wt  T^{-}_{0}}(\wt\om)\in B\,\Big |\,    \ze_{r-1}<\wt  T^{-}_{0}\leq  \ze_{r}\)=1.\label{rkl.17}
 \end{equation}

Since only independent copies of the same process $L^{\cdot}_{i, \ff}$ appear in (\ref{rkl.16}), proofs of the analogues of (\ref{4rev.2}) and (\ref{i.14}), now with $\wt L_{\wt  T_{0}}^{\cdot}$ replaced by $\wt L_{\wt  T^{-}_{0}}^{\cdot}$, are now much simpler.

\section{Behavior of $\wt  L^{x}_{T_{0}} $ at $x$ near $0$ for rebirthed diffusions}\label{sec-atzero}

The following is \cite[Theorem 9.5.25]{book}.

\bt\label{rk-lil}   Let   $\cal Z$ be a recurrent symmetric diffusion in $R^{ 1}$ with
continuous
$\al$-potential densities.    Let $L_{T_0}^x$
denote the local time of $\cal Z$ starting at $y>0$ and killed the first
time it hits
$0$.    Then
\begin{equation}
\limsup_{x\downarrow 0}{  L^{ x}_{T_0}
\over  u_{T_0}(x,x)\log\log (1/u_{T_0}(x,x))}=1
\hspace{.3in}P^y\hspace{.1in}a.s.\label{bm6.1yy}
\end{equation}
  where $u_{T_0}(x,y)=E^x( L^{ y}_{T_0})$.
\et

We want to generalize this to rebirthed processes obtained from diffusions. But note that  the rebirthed process will not be a diffusion since it has jumps at rebirth times.

We assume as in the last section that our renewal measure $\mu$ is supported away from $0$, that is $0\notin \overline {\text{supp}(\mu)}$.   The basic idea is that for a diffusion, the rebirthed process will only be in a neighborhood of $0$ during the rebirth in which $T_0$ occurs, and there we can use (\ref{bm6.1yy}).

\textbf{Example 1:}
 Consider  the symmetric   exponentialy killed  diffusion $Z$ described in Case 4 of \cite[Section 5]{FMR}.  $Z$ is the process with state space $T=R^{1}$ obtained by killing a symmetric recurrent diffusion $\cal Z$   at the end of an independent exponential time with mean $1/\bb$. The potential  density of $Z$ with respect to the speed measure  $m$ is  denoted $\wt u^\bb (x,y)$, see \cite[(5.45)]{FMR}. Let $s(x)$
be the scale function of the diffusion $\cal Z$.  As noted there, the function $u_{T_{0}}(x,y)$ which appears in Theorem \ref{rk-lil} can be expressed in terms of the scale function
\begin{equation}
u_{T_{0}}(x,y)=s(x)\wedge s(y).\label{sc.1}
\end{equation}

Let  $\wt Z\!=\!
(\wt \Om, \wt  \FF_{t}, \wt  Z_t, \wt \th_{t},\wt P^x)$ be the  recurrent Borel right process with state space $ T$ obtained by rebirthing $Z$ as in Section \ref{sec-intro}.   
Let $\wt L_{t}^{x}$ denote the local time of $\wt Z$
  normalized as in (\ref{az.1vw})
and let $L^{x}_{t}  $ be the local time for $Z$ normalized so that
  \begin{equation}
  E^{z}\( L_{\ff}^{x}\)= \wt u^{\bb}(z,x), \quad \forall z,x\in T.
  \end{equation}

\bt\label{theo-nsdiff2} For all $y> 0$,
 \begin{equation}
 \limsup_{x\downarrow 0}{\wt  L^{x}_{\wt T_{0}(\wt \om)}(\wt \om)
\over s(x)\log\log (1/s(x))}=1 \qquad \wt P ^{y}\,\, a.s.\label{85.40}
\end{equation}
 \et

 \textbf{Proof of Theorem \ref{theo-nsdiff2}: } As mentioned after (\ref{74.2}) 
 we cannot have $\wt T_{0}(\wt\om)=\ze_{r}$ for any $r$. Hence by the law of total probability it suffices to show that 
 \begin{equation}
 \wt P ^{y}\( \limsup_{x\downarrow 0}{\wt  L^{x}_{\wt T_{0}(\wt \om)}(\wt \om)
\over s(x)\log\log (1/s(x))}=1 \,\Big |\,    \ze_{r-1}<\wt T_{0}(\wt \om)< \ze_{r} \)=1\label{85.20}
 \end{equation}
 for each $r\geq 1$.
 
 If $\ze (  \om)< T_{0}(  \om)= \inf\{s>0\,|\, \YY_{s}=0 \}$, $Z$ will die at an exponential time before getting close to $0$. Thus, 
with the notation from \cite[Section 2]{FMR}, conditional on $\ze(\om_{i})<T_{0}(\om_{i})$
we must have   for some $\ep(\om_{i})>0$
\be
L^{x}_{\ze(\om_{i})}(\om_{i})=0,\,\hspace{.2in}\forall x\leq \ep(\om_{i}), 
\hspace{.3in}P^y\hspace{.1in}a.s.\label{85.1}
\end{equation}
 Hence by (\ref{74.2}) 
\begin{eqnarray}
&&\hspace{-.5 in}\wt P ^{y}\( \limsup_{x\downarrow 0}{\wt  L^{x}_{\wt T_{0}(\wt \om)}(\wt \om)
\over s(x)\log\log (1/s(x))} =1\,\Big |\,    \ze_{r-1}<\wt T_{0}(\wt \om)<\ze_{r} \)
\nn
\\
&&\hspace{-.5 in}=P ^{\mu}\(\limsup_{x\downarrow 0}{ L^{x}_{T_{0}(\om_{r})}(\om_{r})
\over  s(x)\log\log (1/s(x))} =1\,\Big |\,   T_{0}(\om_{r})<\ze(\om_{r}) \).
\label{85.2}
\end{eqnarray}
(When $r=1$ we have $P ^{y}$ instead of $P ^{\mu}$).
Hence it suffices to show that the last line is $=1$. There we are looking at the total local time up to $T_{0}(\om_{r})$, conditioned 
on not being killed at an exponential time prior to $T_{0}(\om_{r})$. This is given by Theorem \ref{rk-lil} after noting (\ref{sc.1}).\qed

\textbf{Example 2:}
Consider  the symmetric  killed diffusion $Z'$ described in \cite[Section 5, Case 5]{FMR}.  $Z'$ is the  process with state space  $T=(0,\ff) $ that is obtained by starting our symmetric recurrent diffusion $  \cal Z$ in $T$, killing it at the end of an independent exponential time with mean $1/\bb$ and then  killing it  the first time it hits 0. Thus $Z'$ is $Z$ of example 1 killed at  the first time it hits 0. The potential density of $Z'$ with respect to the speed measure  $m$ is denoted $\wt u^{\bb}_{T_{0}}(x,y)$.

Let  $\wt Z\!=\!
(\wt \Om, \wt  \FF_{t}, \wt  Z_t, \wt \th_{t},\wt P^x)$ be the  recurrent Borel right process with state space $ T$ obtained by rebirthing $Z'$ as in Section \ref{sec-intro}.   
Let $\wt L_{t}^{x}$ denote the local time of $\wt Z$
  normalized as in (\ref{az.1vw})
and let $L^{x}_{t}  $ be the local time for $Z'$ normalized so that
  \begin{equation}
  E^{z}\( L_{\ff}^{x}\)= \wt u^{\bb}_{T_{0}}(z,x), \quad \forall z,x\in T.
  \end{equation}
  
     Note that since until its death time, $Z'$ is a diffusion in $T=(0,\ff) $, $\wt Z$ will never hit $0$. As in (\ref{rkl.10}) we set
  \begin{equation}
 \wt T^{-}_{0}( \wt  \om)=\inf \{t\,|\, \wt Z_{t^{-}}(\wt   \om)=0\}. \label{rkl.1}
  \end{equation}

 \bt\label{theo-nsdiff} For all $y>0$,
 \begin{equation}
 \limsup_{x\downarrow 0}{\wt  L^{x}_{ \wt  T^{-}_{0}(\wt \om)}(\wt \om)
\over s(x)\log\log (1/s(x))}=1 \qquad \wt P ^{y}\,\, a.s.\label{85.3}
\end{equation}
 \et

 \textbf{Proof of Theorem \ref{theo-nsdiff}: } By the law of total probability it suffices to show that 
 \begin{equation}
 \wt P ^{y}\( \limsup_{x\downarrow 0}{\wt  L^{x}_{ \wt  T^{-}_{0}(\wt \om)}(\wt \om)
\over  s(x)\log\log (1/s(x))}=1 \,\Big |\,    \ze_{r-1}< \wt  T^{-}_{0}(\wt \om)\leq \ze_{r} \)=1\label{85.20}
 \end{equation}
 for each $r\geq 1$.
 
Since  until its death time, $Z'$ is a diffusion in $T=(0,\ff) $, it will never hit $0$. If $Z'$ dies at an exponential time we must have that for some $\ep(\om_{i})>0$
\be
L^{x}_{\ze(\om_{i})}(\om_{i})=0,\,\hspace{.2in}\forall x\leq \ep(\om_{i}), 
\hspace{.3in}P^y\hspace{.1in}a.s.\label{85.1a}
\end{equation}
Otherwise, $Z'$ dies at the first time $Z$ hits zero. In that case we set 
$ T^{-}_{0}(\om_{i})=    \ze(\om_{i})$.

 Hence by (\ref{74.2r}) 
\begin{eqnarray}
&&\hspace{-.5 in}\wt P ^{y}\( \limsup_{x\downarrow 0}{\wt  L^{x}_{ \wt  T^{-}_{0}(\wt \om)}(\wt \om)
\over  s(x)\log\log (1/s(x))} =1\,\Big |\,    \ze_{r-1}< \wt  T^{-}_{0}(\wt \om)\leq \ze_{r} \)
\nn
\\
&&\hspace{-.5 in}=P ^{\mu}\(\limsup_{x\downarrow 0}{ L^{x}_{T^{-}_{0}(\om_{r})}(\om_{r})
\over  s(x)\log\log (1/s(x))} =1\,\Big |\,  T^{-}_{0}(\om_{r})=\ze(\om_{r})  \).
\label{85.2m}
\end{eqnarray}
Hence it suffices to show that the last line is $=1$. But this follows as in the proof of (\ref{85.2}), since the last line in (\ref{85.2m}) for $Z'$ is exactly the last line in (\ref{85.2}) for $Z$.\qed

\section{ A Generalized Second Ray-Knight Theorem for rebirthed Markov Processes}
\label{sec-NSRK}


Let $ \wt L^{x}_{t}$ be the local time defined in Section \ref{sec-intro} and set \[\wt \tau (t)=\inf\{s>0\,|\,\wt L^{0}_{s}>t \},\] the right continuous inverse of $ \wt L^{0}_{t}$.
In  this section we give a generalized second Ray-Knight Theorem for the local times $\wt L^{x}_{\wt \tau (t)}$ of the non-symmetric Markov process $\wt Z$ which only involves squares of Gaussian processes.
We will also use  
\[\tau_{r} (t)=\inf\{s>0\,|\,L^{0}_{s}(\om_{r})> t \},\] and 
\[\tau^{-}_{r} (t)=\inf\{s>0\,|\,L^{0}_{s}(\om_{r})\geq t \},\]
the right and left continuous inverses of $L^{0}_{s}(\om_{r})$. (We  introduce $\tau^{-}_{r} (t)$  to tie in with the generalized second Ray-Knight Theorem for transient Markov processes, see (\ref{9it4.1t})).

Let $\la,  \la_{1}, \la_{2},\ldots  $ denote  independent exponential random variables of mean $1/p$
 independent of everything else. For a probability measure $P$   we write 
 $P_{\la}=P\times \,d\la$ and similarly, 
 $P_{\la_{i}}=P\times \,d\la_{i}$.
 
 \bl\label{lem-tau}When $\wt L^{0}_{\ze_{r-1}}<t< \wt L^{0}_{\ze_{r}}$ with $r\geq 2$, then 
     \begin{equation}
 \wt L^{x}_{\wt \tau (t)}(\wt\om)=\sum^{r-1}_{i=1}  L^{x}_{\ze(\om_{i})}(\om_{i})+L^{x}_{T_{0}(\om_{r})}(\om_{r})+L^{x}_{\tau_{r} (t- \wt L^{0}_{\ze_{r-1}})
}(\om_{r})\circ \th_{T_{0}(\om_{r})},\label{n94.2fm}
 \end{equation}
 $\wt P^{y}$ a.s. and 
      \begin{equation}
 \wt L^{x}_{\wt \tau (\la)}(\wt\om)=\sum^{r-1}_{i=1}  L^{x}_{\ze(\om_{i})}(\om_{i})+L^{x}_{T_{0}(\om_{r})}(\om_{r})+L^{x}_{\tau^{-}_{r} ((\la- \wt L^{0}_{\ze_{r-1}})
\wedge L^{0}_{\ff}(\om_{r}))}(\om_{r})\circ \th_{T_{0}(\om_{r})},\label{n94.2fn}
 \end{equation}
 for $\wt L^{0}_{\ze_{r-1}}<\la \leq \wt L^{0}_{\ze_{r}}, \quad \wt P_{\la}^{y}$ a.s. 
 \el
 
 \Proof We first note that for any $0\leq a<b$,
  \begin{equation}
\{\wt L^{0}_{a}<t< \wt L^{0}_{b}\}=\{a<\wt \tau (t)<b\}, \label{tau.4}
 \end{equation}
 so that
 \begin{equation}
\{\wt L^{0}_{\ze_{r-1}}<t< \wt L^{0}_{\ze_{r}}\}=\{\ze_{r-1}<\wt \tau (t)< \ze_{r}\}. \label{tau.5}
 \end{equation}
 By (\ref{74.1}), if $\ze_{r-1}<\wt \tau (t)< \ze_{r}$ 
 we have
  \begin{equation}
 \wt L^{x}_{\wt \tau (t)}(\wt\om)=\sum^{r-1}_{i=1}  L^{x}_{\ze(\om_{i})}(\om_{i})+   L^{x}_{\wt \tau (t)-\ze_{r-1}}(\om_{r}),\label{n94.1}
 \end{equation}
 with
 \begin{equation}
0<\wt \tau (t)-\ze_{r-1}< \ze(\om_{r }).\label{tau.6}
 \end{equation}
 We see from this that
 \begin{equation}
\wt \tau (t)-\ze_{r-1}=\tau_{r} (t- \wt L^{0}_{\ze_{r-1}})\label{tau.7}
 \end{equation}
 and consequently using (\ref{tau.7})
  \begin{equation}
0<\tau_{r} (t- \wt L^{0}_{\ze_{r-1}})< \ze(\om_{r }).\label{tau.8}
 \end{equation}
 
 Hence 
 \begin{equation}
L^{x}_{\wt \tau (t)-\ze_{r-1}}(\om_{r})=L^{x}_{T_{0}(\om_{r})}(\om_{r})+L^{x}_{\tau_{r} (t- \wt L^{0}_{\ze_{r-1}})
}(\om_{r})\circ \th_{T_{0}(\om_{r})},\label{tau.9}  
 \end{equation}
 which gives (\ref{n94.2fm}).
 
By  \cite[Lemma 3.6.18]{book}, for each $a>0$
\begin{equation}
\tau_{r} (a)=\tau^{-}_{r} (a), \qquad P_{r}^{0}\quad a.s.\label{tau.10}
\end{equation} 
It follows from this that for each $t>\wt L^{0}_{\ze_{r-1}}$,
\begin{equation}
\tau_{r} (t- \wt L^{0}_{\ze_{r-1}})=\tau^{-}_{r} (t- \wt L^{0}_{\ze_{r-1}}), \qquad P_{r}^{0}\quad a.s.\label{tau.10}
\end{equation}
We obtain from 
(\ref{n94.2fm}) that 
$\wt P_{\la}^{y}$ a.s.,  for $\wt L^{0}_{\ze_{r-1}}<\la <\wt L^{0}_{\ze_{r}}$,
      \begin{equation}
 \wt L^{x}_{\wt \tau (\la)}(\wt\om)=\sum^{r-1}_{i=1}  L^{x}_{\ze(\om_{i})}(\om_{i})+L^{x}_{T_{0}(\om_{r})}(\om_{r})+L^{x}_{\tau^{-}_{r} (\la- \wt L^{0}_{\ze_{r-1}})}\circ \th_{T_{0}(\om_{r})}.
 \label{n94.2fny}
 \end{equation}
Here we used Fubini's Theorem, which also shows that $\wt P_{\la}^{y}\( \la= \wt L^{0}_{\ze_{r}} \)=0$.  But $\wt L^{0}_{\ze_{r-1}}<\la \leq \wt L^{0}_{\ze_{r}}$ shows that $\la -\wt L^{0}_{\ze_{r-1}}\leq L^{0}_{\ff}(\om_{r})$. 
 (\ref{n94.2fn}) now follows.
 \qed

 Set
 \begin{equation}
 F^{x}\(t, \om_{r}\)=L^{x}_{T_{0}(\om_{r})}(\om_{r})+L^{x}_{\tau^{-}_{r} (t
\wedge L^{0}_{\ff}(\om_{r}))}(\om_{r})\circ \th_{T_{0}(\om_{r})}.\label{nF}
 \end{equation}
By (\ref{n94.2fn}),   when $\wt L^{0}_{\ze_{r-1}}<\la\leq \wt L^{0}_{\ze_{r}}$ with $r\geq 2$, 
 we have  
    \begin{equation}
 \wt L^{x}_{\wt \tau (\la)}(\wt\om)=\sum^{r-1}_{i=1}  L^{x}_{\ze(\om_{i})}(\om_{i})+ F^{x}\(\la- \wt L^{0}_{\ze_{r-1}}, \om_{r}\),\quad \wt P_{\la}^{y}\,\,a.s.\label{n94.2f1}
 \end{equation}
  Let
 \bea
 \mathcal{L}_{i}&=&\{L^{x}_{\ze(\om_{i})}(\om_{i}), x\in S\}, \text{ for }1\le i\leq r-1, \text{ and }\label{}\\
 \mathcal{L}'_{r}&=&\nn \{F^{x}\(\la- \wt L^{0}_{\ze_{r-1}}, \om_{r}\), x\in S\},\label{}\\
   \ov {\mathcal{L}} _{r} &=&\nn \{F^{x}\(\la_{r}, \om_{r}\), x\in S\}.\label{}\eea
By (\ref{n94.2f1}) 
     \begin{equation}
 \wt L^{x}_{\wt \tau (\la)}(\wt\om)=\sum^{r-1}_{i=1}   \mathcal{L}_{i}  +  \mathcal{L}'_{r}.\label{n70.2a}
 \end{equation}

 \bl\label{lem-ncondind}
 For all $r\geq 2$ and all Borel sets, $B_1, B_2,\ldots, B_r$, in $C(S, R^1)$,
  \begin{eqnarray}
 &&\wt P_{\la}^{y}\( \mathcal{L}_{1}\in B_1, \ldots,  \mathcal{L}_{r-1}\in B_{r-1},  \mathcal{L}'_{r} \in B_r \,\Big |\,  \wt L^{0}_{\ze_{r-1}}<\la\leq \wt L^{0}_{\ze_{r}}\)
\nn
 \\
 &&\qquad= P_{1,\la_{1}}^{y}\( \mathcal{L}_{1}\in B_1 \,\Big |\,  L^{0}_{\ze (\om_{1})}<\la_{1}   \) \prod^{r-1}_{i=2}   P_{i,\la_{i}}^{\mu}\(\mathcal{L}_{i}\in B_i\,\Big |\,  L^{0}_{\ze (\om_{i})}<\la_{i} \)\nn\\
 &&\qquad\qquad   \times      P_{r,\la_{r}}^{\mu}\( \ov {\mathcal{L}}_{r}\in B_r\,\Big |\,  \la_{r}\leq L^{0}_{\ze (\om_{r})}\).
 \label{n70.5mp}
 \end{eqnarray}
 In particular,
 \begin{equation} \label{n3.26}
   \wt P_{\la}^{y}\(  \mathcal{L} '_{r}\in B_r\,\Big |\,  \wt L^{0}_{\ze_{r-1}}<\la\leq \wt L^{0}_{\ze_{r}}\)=  P_{r,\la_{r}}^{\mu}\( \ov {\mathcal{L}}_{r}\in B_r\,\Big |\,  \la_{r}\leq  L^{0}_{\ze (\om_{r})}\).
\end{equation}
  \el
    (We can state (\ref{n70.5mp}) more efficiently as:  For each $r\geq 2$, $ $
 \begin{equation}
 \mathcal{L}_{1}, \ldots,  \mathcal{L}_{r-1},  \mathcal{L}'_{r}, \label{n70.0}
 \end{equation}
 are `conditionally independent', given $ \wt L^{0}_{\ze_{r-1}}<\la\leq \wt L^{0}_{\ze_{r}}$.)

 \medskip
   \Proof  As in the proof of Lemma \ref{lem-condind0}, let $\mathcal{C}$ be a countable subset of $S$ with compact closure.    For any function $f(x)$,   $x\in\CC$, and $i=1,\ldots, r$,   let  
\begin{equation} \label{}
   (\nu_{i},f)=\sum_{x\in \mathcal{C}}a_{i, x} f(x)
\end{equation}  where the $a_{i, x}\in R^{1}$   and $\sum_{x\in \mathcal{C}}|a_{i, x}|<\ff$  for each $i$. 
Then
 \begin{eqnarray}
 &&\hspace{-.3 in}\wt E_{\la}^{y}\(\exp\({-\sum^{r-1}_{i=1} (\nu_{i},   L^{\cdot}_{\ze(\om_{i})}(\om_{i}))-   \(\nu_{r}, F^{\cdot}\(\la- \wt L^{0}_{\ze_{r-1}}, \om_{r}\)\)}\);    \wt L^{0}_{\ze_{r-1}}<\la\leq \wt L^{0}_{\ze_{r}}\)
\nn
 \\
 &&=\wt E^{y} \int_{\wt L^{0}_{\ze_{r-1}}}^{\wt L^{0}_{\ze_{r-1}}+L^{0}_{\ze (\om_{r})} (\om_{r})} pe^{-pt} \label{n70.3aa}\\
 &&\hspace{.6 in}\exp\({-\sum^{r-1}_{i=1} (\nu_{i},  L^{\cdot}_{\ze(\om_{i})}(\om_{i}))-   \(\nu_{r}, F^{\cdot}\(t- \wt L^{0}_{\ze_{r-1}}, \om_{r}\)\)}     \) \,dt 
 \nonumber \\
 &&=\wt E^{y}\int_{ 0}^{L^{0}_{\ze (\om_{r})}} pe^{-p(u+ \wt L^{0}_{\ze_{r-1}})}\nn\\
 &&\hspace{1 in}  \exp \(-\sum^{r-1}_{i=1} (\nu_{i},  L^{\cdot}_{\ze(\om_{i})(\om_{i}))}-   \(\nu_{r}, F^{\cdot}\(u, \om_{r}\)\)   \)  \,du.\nn    
 \eea
This is equal to,
 \bea
 \lefteqn{\wt E^{y}\( e^{-p \wt L^{0}_{\ze_{r-1}}} e^{-\sum^{r-1}_{i=1} (\nu_{i},  L^{\cdot}_{\ze(\om_{i})}(\om_{i}))}     \)       E^{\mu}\(\int_{ 0}^{L^{0}_{\ze (\om_{r})}} pe^{-pu} e^{-   \(\nu_{r}, F^{\cdot}\(u, \om_{r}\)\) }   \,du\)
 \nonumber}\\
  &&=\wt E^{y}\( \prod^{r-1}_{i=1}  e^{-p L^{0}_{\ze (\om_{i})}} e^{-  (\nu_{i},  L^{\cdot}_{\ze(\om_{i})}(\om_{i}))}     \)  \nn\\
  &&\hspace{1 in}       E^{\mu}\(\int_{ 0}^{L^{0}_{\ze (\om_{r})}} pe^{-pu} e^{- \(\nu_{r}, F^{\cdot}\(u, \om_{r}\)\) }   \,du\)
 \nonumber\\
 &&=  E_{1,\la_{1}}^{y}\(  e^{-  (\nu_{1},  L^{\cdot}_{\ze(\om_{1})}(\om_{1}))} ;  L^{0}_{\ze (\om_{1})}<\la_{1}   \) \prod^{r-1}_{i=2}   E_{i,\la_{i}}^{\mu}\( e^{-  (\nu_{i},  L^{\cdot}_{\ze(\om_{i})}(\om_{i}))} ; L^{0}_{\ze (\om_{i})}<\la_{i} \) \nonumber\\
 &&\hspace{1 in} \times      E_{r,\la_{r}}^{\mu}\( e^{- \(\nu_{r}, F^{\cdot}\(\la_{r}, \om_{r}\)\) }; \la_{r}\leq L^{0}_{\ze (\om_{r})}\).\label{n70.3a}
 \end{eqnarray}
 
By (\ref{n70.3aa})-(\ref{n70.3a}), on setting all $\nu_{i}=0$,
 \bea
 &&\wt P_{\la}^{y} \(\wt L^{0}_{\ze_{r-1}}<\la\leq \wt L^{0}_{\ze_{r}}\)\nn\\
 &&=   P_{1,\la_{1}}^{y} \(L^{0}_{\ze (\om_{1})}<\la_{1} \)\prod^{r-1}_{i=2}   P_{i,\la_{i}}^{\mu} \(L^{0}_{\ze (\om_{i})}<\la_{i}\)   P_{r,\la_{r}}^{\mu} \(  \la_{r}\leq L^{0}_{\ze (\om_{r})}\),\nn
 \eea
 Consequently,
 \begin{eqnarray}
 &&\hspace{-.3 in}\wt E_{\la}^{y}\(\exp\({-\sum^{r-1}_{i=1} (\nu_{i},   L^{\cdot}_{\ze(\om_{i})}(\om_{i}))-   \(\nu_{r}, F^{\cdot}\(\la- \wt L^{0}_{\ze_{r-1}}, \om_{r}\)\)}\)\,\Big |\,     \wt L^{0}_{\ze_{r-1}}<\la\leq \wt L^{0}_{\ze_{r}}\)\nn
 \\
 &&=  E_{1,\la_{1}}^{y}\(  e^{-  (\nu_{1},  L^{\cdot}_{\ze(\om_{1})}(\om_{1}))} \,\Big |\, L^{0}_{\ze (\om_{1})}<\la_{1}   \) \prod^{r-1}_{i=2}   E_{i,\la_{i}}^{\mu}\( e^{-  (\nu_{i},  L^{\cdot}_{\ze(\om_{i})}(\om_{i}))} \,\Big |\, L^{0}_{\ze (\om_{i})}<\la_{i} \) \nonumber\\
 &&\hspace{1 in} \times      E_{r,\la_{r}}^{\mu}\( e^{- \(\nu_{r}, F^{\cdot}\(\la_{r}, \om_{r}\)\) } \,\Big |\, \la_{r}\leq L^{0}_{\ze (\om_{r})}\).
 \nonumber
 \end{eqnarray}
 
 By Laplace inversion this gives  (\ref{n70.5mp}) with the $B_{i}$  restricted to $\mathcal{C}$, and the general case then follows by continuity.
 \qed

 Let 
 $L^{x}_{i, t}$ for $i=1, 2, \ldots$ be independent copies of $L^{x}_{t}$. Set
  \begin{equation}
 F^{x}\(r,s\)=L^{x}_{r,T_{0}} +L^{x}_{r, \tau^{-}_{r} (s
\wedge L^{0}_{r,\ff} )}\circ \th_{T_{0}}.\label{nFr}
 \end{equation}
 
 For any set $C$, let $F(C)$ denote the set of real-valued functions
$f$ on $C$. Define the evaluations
$\it{i}_{ x}:F(C)\mapsto R^{{\rm 1}}$ by $\it{i}_{ x}( f)=f( x)$. We use
$\mathcal{M}( F(C))$ to denote  the smallest $\si$-algebra for which the
evaluations
$\it{i}_{ x}$ are Borel measurable for all
$x\in C$. $\mathcal{M}( F(C))$ is generally referred to  as the
$\si$-algebra  of cylinder sets in $F(C)$. 
 

 \bt\label{theo-n1.1} If for some $r\geq 2$ and measurable set of functions $B\in \mathcal{M}( F(S))$
 \begin{equation}
\(P_{1}^{y}\times \prod^{r-1}_{i=2} P_{i}^{\mu}\times   P_{r, \la_{r}}^{\mu}\)  \( \sum^{r-1}_{i=1} L^{\cdot}_{i, \ff}+ F^{\cdot}\(r, \la_{r}\)\in B\)=1,  \label{n70.6}
 \end{equation}
 then  
  \begin{equation}
\wt P_{\la}^{y} \(\wt L^{\cdot}_{\wt\tau (\la)}(\wt\om)\in B\,\Big |\,   \wt L^{0}_{\ze_{r-1}}<\la\leq \wt L^{0}_{\ze_{r}}\)=1.\label{n70.7}
 \end{equation}
 \et
 
 \noindent\textbf{Proof of Theorem \ref{theo-n1.1} }
 Since $\la_{i}, i=1,\ldots, r-1$ do not appear in the event  $\{ \sum^{r-1}_{i=1} L^{\cdot}_{i, \ff}+ F^{\cdot}\(r, \la_{r}\)\in B\}$, we can just as well write  (\ref{n70.6}) as
  \begin{equation}
\(P_{1,\la_{1}}^{y}\times \prod^{r-1}_{i=2} P_{i,\la_{i}}^{\mu}\times   P_{r, \la_{r}}^{\mu}\)  \( \sum^{r-1}_{i=1} L^{\cdot}_{i, \ff}+ F^{\cdot}\(r, \la_{r}\)\in B\)=1. \label{n70.6c}
 \end{equation}

With this formulation, using the fact that  an almost sure event   with respect to a given probability occurs almost surely for any conditional version of the probability, we see 
  that if (\ref{n70.6c}) holds,   then 
  \be\(P_{1,\la_{1}}^{y}\times \prod^{r-1}_{i=2} P_{i,\la_{i}}^{\mu}\times   P_{r, \la_{r}}^{\mu}\) \( \sum^{r-1}_{i=1} L^{\cdot}_{i, \ff}+F^{\cdot}\(r, \la_{r}\)\in B\,\Big |\,\mathcal{A}\)=1,
  \label{n70.6d}
\ee 
 where
  \begin{equation}
 \mathcal{A}=\{  L^{0}_{1, \ff}<\la_{1},\cdots,   L^{0}_{r-1, \ff}<\la_{r-1},  \,  \la_{r}\leq  L^{0}_{r, \ff}\}.    \label{}
  \end{equation}
 We note that   for each $i$, $\mathcal{L}_{i}$, which is $L^{x}_{\ze(\om_{i})}(\om_{i})$, has the law of $L^{x}_{i,\ff}$, 
 and $ \bar {\mathcal{L}}$, which is $ F^{x}\(\la_{r}, \om_{r}\)$,  has the law of $F^{x}\(r, \la_{r}\)$. Hence (\ref{n70.6d}) is equivalent to 
    \begin{eqnarray}
  &&\(P_{1,\la_{1}}^{y}\times \prod^{r-1}_{i=2} P_{i,\la_{i}}^{\mu}\times   P_{r, \la_{r}}^{\mu}\)
  \label{n70.6e}
  \\
  &&  \(   \sum^{r-1}_{i=1} \mathcal{L}_{i}+\bar {\mathcal{L}}\in B\,\Big |\,    L^{0}_{\ze(\om_{1})}<\la_{1},\cdots,   L^{0}_{\ze(\om_{r-1})}<\la_{r-1},  \,  \la_{r}\leq  L^{0}_{\ze(\om_{r})} \)=1.
  \nonumber
  \end{eqnarray}
It then follows from (\ref{n70.5mp}) that, 
  \begin{equation}
  \wt P_{\la}^{y}\(\sum^{r-1}_{i=1} \mathcal{L}_{i}+\mathcal{L}'_{r}\in B   \,\Big |\,      \wt L^{0}_{\ze_{r-1}}<\la\leq \wt L^{0}_{\ze_{r}} \)=1.\label{}
  \end{equation}
Using (\ref{n70.2a})  we see that this is the statement in (\ref{n70.7}).
\qed

Note that (\ref{n70.6}) will certainly hold if for all $t> 0$
 \begin{equation}
\(P_{1}^{y}\times \prod^{r-1}_{i=2} P_{i}^{\mu}\times   P_{r}^{\mu}\)  \( \sum^{r-1}_{i=1} L^{\cdot}_{i, \ff}+  F^{\cdot}\(r, t\)\in B\)=1.  \label{n70.6x}
 \end{equation}
Furthermore, using 
  \begin{equation}
F^{x}\(r, t\)=L^{x}_{r, T_{0} } +L^{x}_{r,\tau^{-}_{r} (t\wedge   L^{0}_{r, \ff}  )}\circ\th_{T_{0}},\label{n94.2ee}
 \end{equation}
we see that (\ref{n70.6x}) is equivalent to
 \begin{equation}
\(P_{1}^{y}\times \prod^{r}_{i=2} P_{i}^{\mu}\times   P_{r'}^{0}\)  \( \sum^{r-1}_{i=1} L^{\cdot}_{i, \ff}+  L^{\cdot}_{r, T_{0}}+ L^{\cdot}_{r',\tau^{-}_{r'} (t\wedge   L^{0}_{r', \ff}  )}\in B\)=1.  \label{n70.6y}
 \end{equation}
 
  Let $\eta_{i, 0} (x)$, i=1,2,\ldots be independent  Gaussian processes with covariance $ u^{0}(x,y)$ and let $\ov\eta _{i}(x) =1,2$ be independent  Gaussian processes with covariance $  u_{T_{0}}(x,y)$ independent of the $\eta_{i, 0} (x)$. Let $\rho$ be an exponential random variable with mean $u^{0}(0,0)$ that is
independent of everything else and set $h_x=P^{x}(T_{0}<\ff)=\frac{u^{0}(x,0)}{u^{0}(0,0)}$.

Let
\be 
\wh G_{r, s}\(x\)=\sum_{i=1}^{r-1} \frac{1}{2}(\eta_{0,i }(x)+s)^2+\frac{1}{2}(\ov\eta_{1}(x)+s)^2+\frac{1}{2}\ov\eta_{2}^2(x), \label{92.21p}
\ee
\be 
\ov G_{r, s, t}\(x\)=\sum_{i=1}^{r-1} \frac{1}{2}(\eta_{0,i }(x)+s)^2+\frac{1}{2}(\ov\eta_{1}(x)+s)^2+\frac{1}{2}\(\ov\eta_{2}(x)+h_x\sqrt{2( t\wedge \rho)}\)^2, \label{92.21p}
\ee
and $P_{\ov G_{r, s}}=  \prod_{i=1}^{r-1}  P_{\eta_{0,i} } \times  P_{\ov\eta_{1} }\times P_{\ov\eta_{2}}\times P_{\rho}$.

  The following is our generalization of the second Ray-Knight theorem to rebirthed Markov processes.
\bt\label{theo-ITcond2} For all $r\geq 1$, if for some measurable set of functions $B\in \mathcal{M}( F(S))$  
 \begin{equation}
  P_{\ov G_{r, s}}\(\ov G_{r, s,  t}\in B\)=1\quad \text{  for all $ t>0$},  \label{94.6n}
 \end{equation}
 then  for all  $y$
  \begin{equation}
\(\wt P_{\la}^{y}\times P_{\ov G_{r, s}}\) \(\wt L^{\cdot}_{\wt \tau  (\la)}(\wt\om)+ \wh G_{r, s}\in B\,\Big |\,     \wt L^{0}_{\ze_{r-1}}<\la\leq \wt L^{0}_{\ze_{r}}\)=1.\label{94.7n}
 \end{equation}
 \et
 
\Proof  We use the Eisenbaum Isomorphism Theorem, \cite[Theorem 8.1.1]{book}, which gives the relationships,   
\bea\lefteqn{
\Big\{ L^x_{i,\ff}+\frac{1}{2}(\eta_{0,i} (x)+s)^2\,\,;\,x\in S\,,\,P_{i}^y\times
P_{\eta_{0,i}}\Big\}\label{9it1.2}}\\ &&\qquad\stackrel{law}{=}
\Big\{\frac{1}{2}(\eta_{0,i} (x)+s)^2\,\,;\,x\in S\,,\,\(1+\frac{\eta_{0,i} (y)}s\)P_{\eta_{0,i}}\Big\},\nn
\eea
where $P_{\eta_{0, i}}$ is the probability for $\eta_{0, i}$, and as in (\ref{4it1.2a})
\bea\lefteqn{
\Big\{ L^x_{r, T_{0}}+\frac{1}{2}(\ov\eta_{1} (x)+s)^2\,\,;\,x\in S\,,\,   P_{r}^y\times
P_{\ov\eta}\Big\}\label{4it1.2a2}}\\ &&\qquad\stackrel{law}{=}
\Big\{\frac{1}{2}(\ov\eta_{1} (x)+s)^2\,\,;\,x\in S\,,\,(1+\frac{\ov\eta_{1} (y)}s)P_{\ov\eta_{1}}\Big\}\nn
\eea
where $P_{\ov\eta_{1}}$ is the probability for $\ov\eta_{1}$.  We also use the 
  generalized second Ray--Knight theorem, transient case, \cite[(8.54)]{book}, which says that  for any  $  t>0$,
\bea\lefteqn{ \Big\{  L^{x}_{r',\tau^{-}_{r'} (t\wedge   L^{0}_{r', \ff}  )}+\textstyle{ 1\over 2}\ov\eta_{2}^2(x); 
x\in S,P_{r'}^0\times
P_{\ov\eta_{2}}\}\label{9it4.1t}}\\&&\stackrel{law}{=}
\Big\{\textstyle{ 1\over 2}\!\(\ov\eta_{2}(x)+h_x\sqrt{2( t\wedge \rho)}\)^2; x\in
S,P_{\ov\eta_{2}}\times P_{\rho}\Big\}.\nn
\eea

Set $\eta_{0,i} (\mu)=\int \eta_{0,i} (x)\,d\mu(x)$,  $\ov\eta_{1} (\mu)=\int \ov\eta_{1} (x)\,d\mu(x)$.   It follows from the preceding   relationships     that  for each  $r\geq 2$, and $ t>0$
\bea
&&\label{9it1.3}
\Big\{ \sum^{r-1}_{i=1} L^{x}_{i, \ff}+  L^x_{r, T_{0}}+L^{x}_{r',\tau^{-}_{r'} (t\wedge   L^{0}_{r', \ff}  )}\nn\\
&&\hspace{.5 in}+\sum_{i=1}^{r-1}\frac{1}{2}(\eta_{0,i}(x)+s)^2+\frac{1}{2}(\ov\eta_{1}(x)+s)^2+\frac{1}{2}\ov\eta_{2}^2(x)\,;\,\, x\in S,\,\,  \\
&&\hspace{1 in} \(P_{1}^{y}\times \prod^{r}_{i=2} P_{i}^{\mu}\times P_{r'}^{0}\)     \times\prod_{i=1}^{r-1}  P_{\eta_{0,i} } \times   P_{\ov\eta_{1} }  \times   P_{\ov\eta_{2} }\Big\}\nn\\
 &&\stackrel{law}{=}
\Big\{ \sum_{i=1}^{r-1}\frac{1}{2}(\eta_{0,i}(x)+s)^2+\frac{1}{2}(\ov\eta_{1}(x)+s)^2\nn\\
&&\hspace{1 in}+\frac{1}{2}\(\ov\eta_{2}(x)+h_x\sqrt{2( t\wedge \rho)}\)^2\,\,\,\,\,;\,x\in S\,,\nn\\
&& \hspace{-.4 in}\qquad    \,(1+\frac{\eta_{ 1}(y)}s)P_{\eta_{0,1} } \times\prod_{i=2}^{r-1}  (1+\frac{\eta_{0,i}(\mu)}s)P_{\eta_{0,i} }\times   ( (1+\frac{\ov\eta_{1}(\mu)}s)P_{\ov\eta_{1} })\times P_{\ov\eta_{2}}\times P_{\rho} \Big\}.\nn
\eea
The case of $r=1$ is precisely (\ref{9it1.2}).

As in the proof of Theorem \ref{theo-condind0}, our Theorem then follows from Theorem \ref{theo-n1.1} and (\ref{n70.6y}).
 \qed

\subsection{Application: Exact moduli of continuity at $\wt \tau  (t)$ for the local times of rebirthed Markov processes }

   \bt\label{theo-localmod2} Assume that $[0,1]\subseteq S\subseteq R^{1} $.  Let $\eta_{0} (x)$ be a Gaussian process with covariance $ u^{0}(x,y)$ and let $\ov\eta  (x)$ be a Gaussian process with covariance $ u_{T_{0}}(x,y)$. Set,
  \begin{equation} \label{}
  \si_{0}^{2} (x,y)=  u^{0}(x,x)+ u^{0}(y,y)-2 u^{0}(x,y) 
\end{equation}
and,
  \begin{equation} \label{}
  \ov\si^{2} (x,y)= u_{T_{0}}(x,x)+ u_{T_{0}}(y,y)-2  u_{T_{0}}(x,y). 
\end{equation} 
If   for some $d>0$
   \begin{equation}
 \lim_{x,y\da d} \frac{\ov \si^2(x,y)}{\si_{0} ^2(x,y)}=1\label{19.4am},
   \end{equation}
  and   $\{\eta_{0} 
  (x),x\in [0,1] \}$ has a local modulus of continuity,
  \begin{equation} \label{19.7m}
   \limsup_{x\to 0} \frac{|\eta_{0} (x+d)-\eta_{0}  (d)|} {   \phi(x)} =  1    \qquad a.s.
\end{equation}
for some increasing function $\phi$,
 then for any $t>0$
   \be
   \limsup_{x\to 0}\frac{|\wt L_{\wt \tau  (t)}^{x+d}-\wt L_{\wt \tau  (t)}^{d}|}{\phi(x)  }= \(2\wt L_{\wt \tau  (t)}^{d}\)^{1/2},  \quad   \,\,\, \wt P^{y}\,\,
a.s.\label{9rev.2}
\ee
for all $y\in S$.

\et

 \Proof   
 Let $F( \mathcal{C})$  be a set of  functions on a countable dense   subset $\mathcal{C}$ of $[0,1]$ and set,  
\be
B =\Big\{g\in F( \mathcal{C})\,\Big |  \,\limsup_{x\in \mathcal{C},   x\to 0}\frac{|g (x+d)-  g(d)|}{\phi(x)}=  \(2{g(d)}\)^{1/2} \Big\}.\label{9A.local}
\ee
We first show that
\be
  P_{\ov G_{r, s}} \(\ov G_{r, s,  t} \in B\)=1,\label{9az.2p}
\ee 
where
\bea 
&&\ov G_{r, s,  t}\(x\)=\sum_{i=1}^{r-1} \frac{1}{2}(\eta_{0,i }(x)+s)^2+\frac{1}{2}\(\ov\eta_{1}(x)+s\)^2\nn\\
&&\hspace{1 in}+\frac{1}{2}\(\ov\eta_{2}(x)+h_x\sqrt{2(  t\wedge \rho)}\)^2. \label{92.21p}
\eea
We write
\begin{eqnarray}
&& \(\ov\eta_{2}(x+d)+h_{x+d}\sqrt{2(  t\wedge \rho)}\)^2-\(\ov\eta_{2}(d)+h_{d}\sqrt{2(  t\wedge \rho)}\)^2 
\label{92.21p7}
\\
&&= \(\ov\eta_{2}(x+d)+h_{x+d}\sqrt{2(  t\wedge \rho)}\)^2-\(\ov\eta_{2}(x+d)+h_{d}\sqrt{2(  t\wedge \rho)}\)^2 
\nonumber\\
&&+ \(\ov\eta_{2}(x+d)+h_{d}\sqrt{2(  t\wedge \rho)}\)^2-\(\ov\eta_{2}(d)+h_{d}\sqrt{2(  t\wedge \rho)}\)^2 
\nonumber
\end{eqnarray}
 The middle line can be written as  
 \begin{equation}
  2\ov\eta_{2}(x+d)(h_{x+d}-h_{d})\sqrt{2(  t\wedge \rho)} +(h^{2}_{x+d}-h^{2}_{d}) 2(  t\wedge \rho),  \label{}
\end{equation}
and using \cite[(3.13)]{FMR} we see that  
\begin{equation}
\lim_{   x\to 0}\frac{ 2\ov\eta_{2}(x+d)(h_{x+d}-h_{d})\sqrt{2(  t\wedge \rho)} +(h^{2}_{x+d}-h^{2}_{d}) 2(  t\wedge \rho) }{\phi(x)}=0.\label{}
\end{equation}

For the last line of (\ref{92.21p7}) we set $s'=h_{d}\sqrt{2( t\wedge \rho)}$, and to prove 
(\ref{9az.2p}) it suffices to prove that 
\be
  P_{\ov G_{r, s}} \(\wt G_{r, s,  t} \in B\)=1\label{9az.2pp}
\ee 
where now
\be 
\wt G_{r, s,  t}\(x\)=\sum_{i=1}^{r-1} \frac{1}{2}(\eta_{0,i }(x)+s)^2+\frac{1}{2}\(\ov\eta_{1}(x)+s\)^2+\frac{1}{2}\(\ov\eta_{2}(x)+s'\)^2. \label{92.21pp}
\ee
(\ref{9az.2pp})  follows as in the proof of \cite[Corollary 3.1]{FMR}.

It then follows as in the proof of \cite[Theorem  4.1]{FMR} that, for all $y\in S$,
\begin{equation}
\limsup_{ x\in \mathcal{C},   x\to 0}
 \frac{|\wt   L^{x+d}_{\wt \tau  (\la)}- \wt  L^{ d}_{\wt \tau  (\la)}|}{\phi(x)}  
=  \( 2\wt   L^{ d}_{\wt \tau  (\la)}\)^{1/2},\quad \wt P_{\la}^{y}\quad a.s. \label{}
\end{equation}

 Since  $\wt  L^{x}_t$ is continuous in $x$ we can remove the condition that $x\in \mathcal{C}$. This gives us  (\ref{9rev.2}) for a.e. $t$, $\wt P^{y}$\,\,
a.s. But $ \wt   L^{x+d}_{\wt \tau  (t)}$ is flat in $t$ for $d>0$ and $x$ near $0$, so (\ref{9rev.2}) actually holds for all $t$, $\wt P^{y}$\,\,
a.s.
\qed


\bt\label{theo-unifmod2}    
Assume that $[0,1]\subseteq S\subseteq R^{1} $. With the notation of the last Theorem assume that
  \begin{equation} \label{fi.12}
 \lim_{x\to 0}\sup_{|u-v|\le x}(\si_{0} )^{2} (u,v)\log 1/ |u-v|=0, 
\end{equation}
and similarly for $(\ov \si_{0} )^{2}$.

  Let $\vf(u,v)$ be   such that for some constant $0<C<\ff$,
\begin{equation} \label{fi.13}
  \vf(u,v)\le  \wt\vf\(C|u-v|\),
\end{equation}
 where $\wt\vf$ is an increasing continuous function with $\wt\vf(0)=0$.     If
 \be
 \lim_{h\to 0}\sup_{\stackrel{|u-v|\le h }{ u,v\in\De}} \frac{ ( \eta _{0 }(u)-\eta _{0}(v)) a_{0} + (\ov\eta (u)-\ov\eta (v)) a_{1} }{    \    \vf (u,v) } =1 \label{fi.11} ,
\ee
 for all   intervals $\De$ in $[0,1]$ and all     $\{a_0,a_1\}$  with  $a_0^2+a_1^2=1$     then,  
   \be
\lim_{h\to 0}\sup_{\stackrel{|u-v|\le h }{ u,v\in\De}} \frac{|\wt L_{\wt \tau  (t)}^{u}-\wt L_{\wt \tau  (t)}^{v}|}{\vf (|u-v|)  }= \sup_{u\in \De}\(2\wt L_{\wt \tau  (t)}^{u}\)^{1/2},  \quad    a.e.\,\,\, t,  \,\,\,  \wt P^{y}\,\,
a.s.\label{9urev.2}
\ee
for all $y\in S$. 

If $\De=[c,d]$ with $c>0$ then
   \be
\lim_{h\to 0}\sup_{\stackrel{|u-v|\le h }{ u,v\in\De}} \frac{|\wt L_{\wt \tau  (t)}^{u}-\wt L_{\wt \tau  (t)}^{v}|}{\vf(|u-v|)  }= \sup_{u\in \De}\(2\wt L_{\wt \tau  (t)}^{u}\)^{1/2},  \quad      \,\,\,  \wt P^{y}\,\,
a.s.\label{9urev.3}
\ee
\et

The proof is almost the same as that of Theorem \ref{theo-localmod2} except that instead of $B$ in (\ref{9A.local}) we now use 
\be
A =\Big\{g\in F( \mathcal{C})\,\Big |  \,\lim_{h\to 0}\sup_{\stackrel{|u-v|\le h }{ u,v\in\De\cap \mathcal{C}}}\frac{|g (u)-  g(v)|}{\phi(|u-v|)}= \sup_{u\in \De} \(2{g(u)}\)^{1/2} \Big\},\label{9B.unif}
\ee
\cite[Theorem  3.2]{FMR},
and  the fact that the random variables $\sup_{u\in \De}|\eta_{i }(u)|$, $i=1,\ldots,r-1$,  $\sup_{u\in \De}|\ov \eta_{i} (u)|$, $i=1,2$,  take values arbitrarily close to 0 with probability greater than 0.\qed

{\bf  Acknowledgements: } We want to thank Michael B. Marcus for many helpful suggestions.

       \bibliographystyle{amsplain}

\bigskip
\noindent
\begin{tabular}{lll} &   P.J. Fitzsimmons  \\
&Department of
Mathematics  \\
& University of California,   San Diego\\
&La Jolla CA, 92093 USA  \\ 
& pfitzsim@ucsd.edu\\ 
& &\\
& & \\
& Jay Rosen\\
& Department of Mathematics\\
&  College of Staten Island, CUNY\\
& Staten Island, NY 10314, USA \\
& jrosen30@optimum.net
\end{tabular}

\end{document}